\documentclass[12pt,a4paper,draft]{amsart}
\usepackage[english]{babel}
\usepackage{a4wide,amssymb,xypic,ifthen}

\newcommand{\qee}{ \hfill\hspace{2pt}$\triangle$}
\newcommand{\marginnote}[1]{\ifthenelse{\isodd{\thepage}}{\normalmarginpar}
{\reversemarginpar}\marginpar{\fbox{\parbox{15mm}{\sloppy\footnotesize #1}}}}

\newtheorem{thm}{Theorem}[section]
\newtheorem{corol}[thm]{Corollary}
\newtheorem{lemma}[thm]{Lemma}
\newtheorem{prop}[thm]{Proposition}
\newtheorem{defin}[thm]{Definition}

\theoremstyle{remark}
\newtheorem{rema}[thm]{Remark}
 \newenvironment{remark}{\begin{rema}}{\qee\end{rema}}
\newtheorem{exe}[thm]{Example}

\newcommand{\PP}{{\mathbb P}}

\newcommand{\R}{\mathbb R}
\newcommand{\C}{\mathbb C}
\newcommand{\Q}{\mathbb Q}

\def\ker{\operatorname{ker}}

\def\rk{\operatorname{rk}}
\def\Ad{\operatorname{Ad}}

\def\End{\operatorname{End}}
\def\dim{\operatorname{dim}}
\newcommand\grass{\mbox{Gr}}
\newcommand\hgrass{{\mathfrak{Gr}}}

\newcommand{\fE}{{\mathfrak E}}

\newcommand{\fR}{{\mathfrak R}}

\newcommand{\g}{{\mathfrak g}}

\newcommand{\cat}[1]{\boldsymbol{\operatorname{#1}}}

\begin{document}
\bigskip\bigskip
\title[Semistable   principal Higgs bundles]{Semistable  principal Higgs bundles}
\bigskip
\date{\today}
\thanks{ This research was partly supported by the Spanish {\scriptsize
DGES} through the research project BFM2006-04779 by ``Junta de Castilla y Le{\'o}n''
through the research project SA001A07, by Istituto Nazionale per l'Alta Matematica, by
the Italian National Research Project ``Geometria delle variet\`a algebriche'',   by the European project {\sc misgam} and by
a grant of the University of Salamanca.  Both authors are members of the research
group {\scriptsize VBAC} (Vector Bundles on Algebraic Curves).}

 \maketitle \thispagestyle{empty}
\begin{center}{\sc Ugo Bruzzo} \\
Scuola Internazionale Superiore di Studi Avanzati,\\ Via Beirut 2-4, 34013
Trieste, Italia; \\ Istituto Nazionale di Fisica Nucleare, Sezione di Trieste \\ E-mail: {\tt bruzzo@sissa.it} \\[10pt]
{\sc Beatriz Gra\~na Otero} \\
Departamento de Matem\'aticas  and Instituto de F\'\i sica \\ Fundamental y Matem\'aticas, Universidad de Salamanca,
\\ Plaza de la Merced 1-4, 37008 Salamanca, Espa\~na\\ E-mail: {\tt beagra@usal.es}
\end{center}

\bigskip\bigskip

\begin{abstract} We give a Miyaoka-type semistability criterion for principal Higgs $G$-bundles $\fE$
on complex projective manifolds of any dimension, i.e., we prove that $\fE$ is semistable and the second
Chern class of its adjoint bundle vanishes if and only if certain line bundles, obtained from the characters
of the parabolic subgroups of $G$, are numerically effective. We also give alternative characterizations
in terms of a notion of numerical effectiveness of Higgs vector bundles we have recently introduced. \end{abstract}

\newpage

\section{Introduction.}
In 1987 Miyaoka gave a criterion for the semistability of a vector bundle $V$ on a projective curve in terms
of the numerical effectiveness of a suitable divisorial class (the relative anticanonical divisor of the
projectivization $\mathbb PV$ of $V$). Recently several generalizations of this criterion have been formulated \cite{BH,BB,BSch},
dealing with principal bundles, higher dimensional varieties, and considering also the case of bundles
on compact K\"ahler manifolds. In this paper we prove a Miyaoka-type criterion for principal Higgs bundles on complex
projective manifolds. Let us give a rough anticipation of this result. Given a principal Higgs $G$-bundle $E$
on a complex projective manifold $X$, with Higgs field $\phi$,
 and a parabolic subgroup $P$ of $G$, we introduce a subscheme
${\fR}_P(E,\phi)$ of the total space of the bundle $E/P$ whose sections parametrize reductions of the structure group $G$ to $P$
that are compatible with the Higgs field $\phi$. Then in Theorem \ref{Miyahigher} we prove the equivalence
of the following conditions: for every reduction of $G$ to a parabolic subgroup $P$ which is compatible with the Higgs field, and every dominant character of $P$, a certain associated line bundle on ${\fR}_P(E,\phi)$ is numerically effective; $(E,\phi)$
is semistable as a principal Higgs bundle, and the second Chern class of the adjoint bundle
$\Ad(E)$ (with real coefficients) vanishes. We first prove this fact when $X$ is a curve (so that the
condition involving the second Chern class is void) and then extend it to complex projective manifolds
of arbitrary dimension.

In Section \ref{another}  we   formulate an additional equivalent criterion which states that $\fE=(E,\phi)$ is semistable, and $c_2(\Ad(E))=0$, if and
only if the adjoint Higgs bundle $\Ad(\fE)$ is numerically effective (as a Higgs bundle) in a sense that
we introduced in a previous paper \cite{BG1}.


The final Section \ref{tannaka} develops some Tannakian considerations; basically we show
the equivalence of proving our main theorem \ref{Miyahigher} for principal Higgs bundles
or for  Higgs vector bundles.

Since a principal Higgs bundle with zero Higgs field is exactly a
principal bundle, all results we prove in this paper hold true for
principal bundles. In this way we mostly recover
well-known results or some of the results in \cite{BB} with their
proofs, at other times we provide simpler demostrations, while at
times the results are altogether new.

{\bf Acknowledgement.} This paper has been finalized during a
visit of both authors at the University of Pennsylvania. We thank
Upenn for hospitality and support, and the staff and the
scientists at the Department of Mathematics for providing an
enjoyable and productive atmosphere. We thank M.S.~Narasimhan and Tony Pantev for
valuable suggestions.

\bigskip\section{Semistable principal bundles}\label{Sec:2}
In this short section we recall some basics about principal
bundles, notably the definition of semistable principal bundle (a
basic reference about this topic is \cite{BalSes}). Let $X$ be a
smooth complex projective variety, $G$ a complex reductive
algebraic group, and $\pi\colon E\to X$ a principal $G$-bundle on
$X$. If $\rho\colon G \to \operatorname{Aut}(Y)$ is a
representation of $G$ as automorphisms of a variety $Y$, we may
construct the associated bundle $E(\rho)=E\times_\rho Y$, the
quotient of $E\times Y$ under the   action of $G$ given by
$(u,y)\mapsto (ug,\rho(g^{-1})y)$ for $g\in G$. If $Y=\g$ is the
Lie algebra of $G$, and $\rho$ is the adjoint action of $G$ on
$\g$, one gets the adjoint bundle of $E$, denoted by  $\Ad(E)$.
Another important example is obtained when $\rho$ is given by a
group homomorphism $\lambda\colon G\to G'$; in this case the
associated bundle $E'=E\times_\lambda G'$ is a principal
$G'$-bundle. We say that the structure group $G$ of $E$ has been
extended to $G'$.

If $E$ is a principal $G$-bundle on $X$, and $F$ a principal
$G'$-bundle on $X$, a morphism $E\to F$ is a pair $(f,f')$, where
$f'\colon G\to G'$ is a group homomorphism, and $f\colon E\to F$
is a morphism of bundles on $X$ which is $f'$-equivariant, i.e.,
$f(ug)=f(u)f'(g)$. Note that this induces a vector bundle morphism
$\tilde f\colon \Ad(E)\to\Ad(F)$ given by $\tilde
f(u,\alpha)=(f(u),f'_\ast(\alpha))$, where $f'_\ast\colon\g\to\g'$
is the morphism induced on the Lie algebras.

If $K$ is a closed subgroup of $G$, a \emph{reduction} of the structure
group $G$ of $E$ to $K$ is a principal $K$-bundle $F$ over $X$ together
with an injective $K$-equivariant bundle morphism $F \to E$.
Let $E(G/K)$ denote
the bundle over $X$ with standard fibre $G/K$ associated to $E$ via
the natural action of $G$ on the homogeneous space $G/K$.
There is an isomorphism   $E(G/K)\simeq E/K$ of bundles over $X$. Moreover,
the reductions of the structure group of $E$ to $K$ are in a one-to-one correspondence
with sections $\sigma\colon X \to E(G/K)\simeq E/K$.

 We first recall the definition of
semistable principal bundle when the base variety $X$ is a curve.
Let $T_{E/K,X}$ be the vertical tangent bundle to the bundle
$\pi_K\colon E/K \to X$.

\begin{defin} Let $E$ be a principal $G$-bundle on  a smooth connected projective curve $X$. We say that
$E$ is semistable if for every proper parabolic subgroup $P\subset G$, and every reduction $\sigma\colon X \to E/P$,
the pullback $\sigma^\ast (T_{E/P,X})$ has nonnegative degree.
\end{defin}

When $X$ is a higher dimensional variety, the definition must be somewhat
refined; the introduction of an open dense subset whose complement
has codimension at least two should be compared with the definition
of semistable vector bundle, which involves non-locally free subsheaves
(which are subbundles exactly on open subsets of this kind).

\begin{defin} Let $X$ be a polarized smooth projective variety. A principal  $G$-bundle $E$  on $X$ is semistable
 if and only   if  for any proper parabolic subgroup $P\subset G$, any open dense subset $U\subset X$ such that $\operatorname{codim}(X-U)\ge 2$, and any  reduction $\sigma\colon U \to (E/P)_{\vert U}$ of $G$ to $P$ on $U$, one has $\deg \sigma^\ast (T_{E/P,X}) \ge 0$.
\end{defin}

Here it is important that the smoothness of $X$ guarantees that a
line bundle defined on an open dense subset of $X$ extends
uniquely to the whole of $X$, so that we may consistently consider
its degree. This is discussed in detail in \cite{RamaSub}.

\bigskip\section{Principal Higgs bundles}

We switch now to principal Higgs bundles.
Let $X$ be a smooth complex projective variety, and $G$ a reductive
complex  algebraic group.

\begin{defin} A principal Higgs $G$-bundle $\fE$ is a pair $(E,\phi)$, where
$E$ is a principal $G$-bundle, and $\phi$  is a global section
 of $\Ad(E)\otimes\Omega^1_X$ such that $[\phi,\phi]=0$.
\end{defin}

When $G$ is the general linear group, under the identification $\Ad(E)\simeq
\End(V)$, where $V$ is the vector bundle corresponding to $E$, this agrees with the usual definition  of Higgs vector bundle.

\begin{defin} \label{trivial} A principal Higgs $G$-bundle $\fE=(E,\phi)$ is \emph{trivial}
if $E$ is trivial, and $\phi=0$.
\end{defin}

A morphism between two principal Higgs bundles $\fE=(E,\phi)$ and $\fE'=(E',\phi')$
is a principal bundle morphism $f\colon E\to E'$ such that $(f_\ast\times \operatorname{id})(\phi)=\phi'$, where
$f_\ast\colon \Ad(E)\to\Ad(E')$ is the induced morphism between the adjoint bundles.

Let $K$ be a closed subgroup of $G$, and $\sigma\colon X \to
E(G/K)\simeq E/K$ a reduction of the structure group of $E$ to
$K$. So one has a principal $K$-bundle $F_\sigma$ on $X$ and a
principal bundle morphism $i_\sigma\colon F_\sigma\to E$ inducing
an injective morphism of bundles $\Ad(F_\sigma) \to \Ad(E)$. Let
$\Pi_\sigma\colon  \Ad(E)\otimes
\Omega^1_X\to(\Ad(E)/\Ad(F_\sigma))\otimes \Omega^1_X $ be the
induced projection.

\begin{defin} A section $\sigma\colon X\to E/K $ is a {\em Higgs reduction} of $(E,\phi)$
if $\phi\in\ker \Pi_\sigma$.
\end{defin}

When this happens, the reduced bundle $F_\sigma$ is equipped with a Higgs field $\phi_\sigma$ compatible with $\phi$ (i.e., $(F_\sigma,\phi_\sigma)\to (E,\phi)$ is a morphism of principal Higgs bundles).

\begin{remark} Let us again consider the case when $G$ is the general linear group $Gl(n,\C)$, and let us assume that
$K$ is a (parabolic) subgroup such that $G/K$ is the Grassmann variety $\operatorname{Gr}_k(\C^n)$ of $k$-dimensional quotients
of $\C^n$. If $V$ is the vector bundle corresponding to $E$, a reduction $\sigma$ of $G$ to $K$ corresponds to a rank $n-k$
subbundle $W$ of $V$, and the fact that $\sigma$ is a Higgs reduction means that $W$ is $\phi$-invariant, i.e., $\phi(W)\subset W\otimes \Omega^1_X$.
\end{remark}

We want to show that the choice of $\phi$ singles out a subscheme
of the variety $E/K$, which describes the Higgs reductions of the
pair $(E,\phi)$. Let $\tilde E$ denote the principal $K$-bundle
$E\to E/K$. Since the vertical tangent bundle $T_{E/K,X}$ is the
bundle associated to $\tilde E$ via the adjoint action of $K$ on
the quotient $\mathfrak g/\mathfrak k$, and $\pi_K^\ast \Ad(E)$ is
the bundle associated to $\tilde E$ via the adjoint action of $K$
on $\mathfrak g$, there there is a natural morphism $\chi\colon
\pi_K^\ast \Ad(E) \to T_{E/K,X}$. Then $\phi$ determines a section
$\chi(\phi) := (\chi \otimes \text{id})(\pi_{K}^{*}\phi)$ of
$T_{E/K,X}\otimes \Omega^1_{E/K}$.

\begin{defin} The {\em scheme of  Higgs reductions of $\fE=(E,\phi)$ to $K$}
is the closed subscheme  ${\fR}_K(\fE)$ of $E/K$ given by the zero locus of $\chi(\phi)$.
\end{defin}

This construction is compatible with base change, i.e., if
$f\colon Y \to X$ is a morphism of smooth complex projective
varieties, and $f^\ast(\fE)$ is the pullback of $\fE$ to $Y$, then
${\fR}_K(f^\ast(\fE))\simeq Y\times_X {\fR}_K(\fE)$. By
construction, $\sigma\colon X \to E(G/K)\simeq E/K$ is a Higgs
reduction if and only if it takes values in the subscheme
${\fR}_K(\fE)\subset E/K$. Moreover the scheme of Higgs reductions
is compatible with morphisms of principal Higgs bundles.  This
means that if $\fE=(E,\phi)$ is a principal Higgs $G$-bundle,
$\fE'=(E',\phi')$ a principal Higgs $G'$-bundle, $\psi\colon G\to
G'$ is a group homomorphism, and $f\colon\fE\to\fE'$ is a
$\psi$-equivariant morphism of principal Higgs bundles, then for
every closed subgroup $K\subset G$  the induced morphism $E/K\to
E'/K'$, where $K'=\psi(K)$, maps ${\fR}_K(\fE)$ into
${\fR}_{K'}(\fE')$.

Also, one should note that the scheme of Higgs reductions is in general singular,
so that in order to consider  Higgs bundles on it one needs to use the theory
of the de Rham complex for arbitrary schemes, as developed by Grothendieck
\cite{EGAIV}.

For the time being we restrict our attention to the case when $X$
is a curve. We start by introducing a notion of semistability for
principal Higgs bundles (which is equivalent to the one given in
Definition 4.6 in \cite{AnBis}).

\begin{defin}\label{Hstabcurv} Let $X$ be a smooth projective curve. A principal Higgs $G$-bundle
$\fE=(E,\phi)$ is stable (semistable)  if for every parabolic subgroup $P\subset G$
and every Higgs reduction $\sigma\colon X\to {\fR}_P(\fE)$ one has
$\deg \sigma^\ast (T_{E/P,X})> 0$ ($\deg \sigma^\ast (T_{E/P,X})\ge 0$). \end{defin}

\begin{remark} \label{remSimp} 
A notion of semistability for principal Higgs bundles
was introduced by Simpson in \cite{Si1}. According to that definition,
a principal Higgs $G$-bundle $\fE$ is semistable if there exists
a faithful linear representation $\rho\colon G\to\operatorname{Aut}(W)$
such that that associated Higgs vector bundle $\mathfrak W = \fE\times_\rho W$
is semistable.  This is not equivalent to our definition,   even in the case
of ordinary (non-Higgs) principal bundles (in which case of course our definition coincides with Ramanathan's classical definition of stability for principal bundles \cite{Rama}). Indeed, according to the definition we adopt,
if $T$ is a torus in $Gl(n,\C)$, any principal $T$-bundle is stable. However the vector bundle
associated to it by the natural inclusion $T\hookrightarrow Gl(n,\C)$ (a direct sum of line bundles)
may fail to be semistable. (Note indeed that this inclusion, regarded as a linear representation
of $T$, does not satisfy the condition in part (ii) of Lemma \ref{adjsemi}.) On the other hand,
this is compatible with the Hitchin-Kobayashi correspondence for principal bundles, which states
that a principal $G$-bundle $E$, where $G$ is a connected reductive complex group, is stable
if and only if it admits a reduction of the structure group to the maximal compact subgroup $K$
of $G$ such that the curvature of the  unique connection on $E$ compatible with the reduction takes
values in the centre of the Lie algebra of $K$ \cite{RamaSub}.
\end{remark}

\begin{lemma} Let $f\colon X' \to X$ be a nonconstant morphism
of smooth projective curves, and $\fE$ a principal Higgs $G$-bundle on $X$.
The pullback Higgs bundle $f^\ast\fE$ is semistable if and only if
$\fE$ is.
\label{basechange}
\end{lemma}
\begin{proof} As we shall prove in Lemma \ref{adjsemi}
in the case of $X$ of arbitrary
dimension, a principal Higgs bundle $\fE$ is semistable if and only if the
adjoint Higgs bundle $\Ad(\fE)$ is semistable (as a Higgs vector bundle).
In view of this result, our claim reduces to the analogous statement
for Higgs vector bundles, which was proved in \cite{BH}.
\end{proof}

If $\fE=(E,\phi)$ is a principal Higgs $G$-bundle on $X$, and
$K$ is a closed subgroup of $G$, we may associate with every character $\chi$ of $K$
a line bundle $L_\chi=E\times_\chi \C$ on $E/K$, where we regard $E$ as a principal
$K$-bundle on $E/K$. An elegant way to state results about reductions is to introduce the notion
of \emph{slope} of a reduction: we call $\mu_\sigma$, the slope of a  Higgs reduction $\sigma$,
the group homomorphism $\mu_\sigma\colon \mathcal X(K) \to \Q$ (where $ \mathcal X(K) $
is the group of characters of $K$) which to any character $\chi$ associates the degree
of the line bundle $\sigma^\ast(L_\chi^\ast)$.

The following result extends to Higgs bundles Lemma 2.1 of \cite{Rama},
and its proof follows easily from that of that Lemma.

\begin{prop} A principal Higgs $G$-bundle
$\fE=(E,\phi)$ is semistable  if and only if for every parabolic subgroup $P\subset G$,
every nontrivial dominant character $\chi$ of $P$, and every Higgs reduction $\sigma\colon X\to {\fR}_P(\fE)$, one has $\mu_\sigma(\chi)\ge 0$.
\label{RamaHiggs}\end{prop}
\begin{proof} Let $\alpha_1,\dots,\alpha_r$ be simple roots
of the Lie algebra $\g$. We may assume that $P$ is a maximal parabolic
subgroup corresponding to a root $\alpha_i$. It has been proven in
\cite[Lemma 2.1]{Rama} that the determinant of the vertical tangent bundle
$T_{E/P,X}$ is associated to the principal $P$-bundle $ E \to E/P$
via a character that may be expressed as $\mu = - m\lambda_i$, where
$\lambda_i$ is the weight corresponding to $\alpha_i$, and $m\ge 0$.
Thus, if $\sigma \colon X \to {\fR}_P(\fE)$ is a Higgs reduction,
 $\deg(\sigma^\ast(L^\ast_\mu))\ge 0 $ if and only if $\deg \sigma^\ast (T_{E/P,X})\ge 0$.
 \end{proof}

 \begin{remark} By ``root'' of $\g$ we mean a root of the semisimple part of $\g$ extended by zero on the centre.  \end{remark}

We may now state and prove a Miyaoka-type semistability criterion
for principal Higgs bundles (over projective curves). This generalizes Proposition 2.1
of \cite{BB}, and, of course, Miyaoka's original criterion in \cite{Mi}.

\begin{thm} A principal Higgs $G$-bundle
$\fE=(E,\phi)$ on a smooth projective curve $X$ is semistable if and only if for every parabolic subgroup $P\subset G$, and
every nontrivial dominant character $\chi$ of $P$, the line bundle $L^\ast_\chi$ restricted to
${\fR}_P(\fE)$ is nef. \label{princMiya}
\end{thm}
\begin{proof} Assume that $\fE$ is semistable and that ${L^\ast_\chi}_{\vert {\fR}_P(\fE)}$
is not nef. Then there is an irreducible curve $Y\subset {\fR}_P(\fE)$ such that
$[Y]\cdot c_1(L^\ast_\chi)<0$. Since $\chi$ is dominant, the line bundle $L^\ast_\chi$ is nef when restricted
to a fibre of the projection $E/P\to X$, so that the curve $Y$ cannot be contained in such a fibre.
Then $Y$ surjects onto $X$. One can choose  a morphism of  smooth
projective curves  $h\colon Y'\to X$   such that $\tilde Y = Y' \times_X Y$ is a disjoint union of smooth curves in $h^\ast({\fR}_P(\fE))$, each mapping isomorphically to $X$.
Using Lemma \ref{basechange}, in this way we may assume that $Y$ is the image of a section
$\sigma \colon X \to {\fR}_P(\fE)$.   Proposition
\ref{RamaHiggs} imples that $[Y]\cdot c_1(L^\ast_\chi)\ge 0$, but this contradicts our assumption.

The converse is obvious in view of Proposition \ref{RamaHiggs}.
\end{proof}

\begin{remark} \label{remframe} Let  $G$ be the linear group $\operatorname{Gl}(n,\C)$.
 If $\fE=(E,\phi)$ is a principal
Higgs $G$-bundle, and $V$ is  the rank $n$ vector bundle corresponding to $E$,
then the identification $\Ad(E)\simeq\End(V)$ makes $\phi$ into a Higgs
morphism $\tilde\phi$  for $V$. A simple calculation shows that the semistability
of $\fE$ is equivalent to the semistability of the Higgs vector bundle
$(V,\tilde\phi)$.

If $P$ is such that the quotient $G/P$ is the $(n-1)$-dimensional
projective space, the bundle  $E/P$ is isomorphic to the projectivization
$\PP V\to X$ of $V$ (regarded as the space whose sections classify
rank 1 locally-free quotients of $V$).
More generally, let $P_k$ be a maximal parabolic subgroup, so that
$G/P_k$ is the Grassmannian of rank $k$  quotient spaces of $\C^n$ for some $k$. In this case
$E/P_k$ is the Grassmann bundle $\operatorname{Gr}_k(V)$ of rank $k$ locally free
quotients of $V$.
Then Theorem \ref{princMiya} corresponds to the result given in \cite{BH},
according to which $(V,\phi)$ is semistable if and only if certain numerical classes
$\theta_k$ in a closed subscheme of $\operatorname{Gr}_k(V)$ are nef (see  \cite{BH,BG1,BG2} for details).
\end{remark}

\bigskip\section{The higher-dimensional case}
In this section we consider the case of a base variety $X$ which is
a complex projective manifold of any dimension. Let $X$ be equipped with a polarization $H$, and let $G$ be a reductive
complex algebraic group.

\begin{defin} \label{Higgsstabhigher} A principal Higgs $G$-bundle $\fE=(E,\phi)$ is stable (semistable)
if and only if for any proper parabolic subgroup $P\subset G$,   any open dense subset $U\subset X$ such that
$\operatorname{codim}(X-U)\ge 2$, and any Higgs reduction $\sigma\colon U \to {{\fR}_P(\fE)}_{\vert U}$ of $G$ to $P$ on $U$,
one has $\deg \sigma^\ast (T_{E/P,X}) \ge 0$. \end{defin}

\begin{remark} The arguments in the proof of Proposition \ref{RamaHiggs}
go through also in the higher dimensional case, allowing one to show that
a principal Higgs $G$-bundle $\fE$ is semistable (according to Definition
\ref{Higgsstabhigher}) if and only if for any proper parabolic subgroup $P\subset G$, any nontrivial dominant character $\chi$ of $P$, any open dense subset $U\subset X$ such that
$\operatorname{codim}(X-U)\ge 2$, and any Higgs reduction $\sigma\colon U \to {{\fR}_P(\fE)}_{\vert U}$ of $G$ to $P$ on $U$,
the line bundle $\sigma^\ast( L_{\chi}^\ast)$ has positive (nonnegative) degree.
\label{remRama}
 \end{remark}

If $\fE$ is a principal Higgs $G$-bundle, we denote by
$\Ad(\fE)$ the Higgs vector bundle given by the adjoint bundle
$\Ad(E)$ equipped with the induced Higgs morphism.

We also introduce the notion of extension of the structure group
for a principal Higgs $G$-bundle $\fE=(E,\phi)$. Given a group
homomorphism $\lambda\colon G\to G'$, we consider the extended
principal bundle $E'$. The group $G$ acts on the Lie algebra $\g'$
of $G'$ via the homomorphism $\lambda$ (and the adjoint action of
$G$), and the $\g'$-bundle associated to $E$ via the adjoint
action of $G'$ is isomorphic to $\Ad(E')$. In this way the Higgs
field of $\fE$ induces a Higgs field for $\fE'$. More generally,
if $\rho\colon G\to\operatorname{Aut}(V)$ is a linear
representation of $G$, the Higgs field of $\fE$ induces a Higgs
field for the associated vector bundle $E\times_\rho V$.

It is known that certain  extensions of the structure group of a
semistable principal bundle are still semistable \cite{RamaRama},
and that a principal bundle is semistable if and only if its
adjoint bundle is \cite{Rama}. The same is true in the Higgs case.

\begin{lemma} \label{adjsemi} (i) A principal Higgs bundle
$\fE$ is semistable if and only if $\Ad(\fE)$ is semistable (as a Higgs vector bundle).

(ii) A principal Higgs $G$-bundle $\fE=(E,\phi)$ is semistable if and only if
for every  linear representation $\rho\colon G\to\operatorname{Aut}(V)$ of $G$ such that
$\rho(Z(G)_0)$ is contained in the centre of $\operatorname{Aut}(V)$, the associated Higgs vector
bundle $\mathfrak V = \fE\times_\rho V$ is semistable (here $Z(G)_0$ is the
component of the centre of $G$ containing the identity).
\end{lemma}
\begin{remark} Let us at first note that if $G$ is the general linear group $Gl(n,\C)$,
the first claim  holds true quite trivially: $\fE$ is semistable if and only if the corresponding
Higgs vector bundle $\mathfrak V$ is semistable, and one knows that $\Ad(\fE)\simeq\End(\mathfrak V)$
is semistable if and only if $\mathfrak V$ is.\end{remark}
\begin{proof}
The first claim is Lemma 4.7 of \cite{AnBis}. The second claim is proved as in
Lemma 1.3 of \cite{AAB}.
\end{proof}

\begin{prop} \label{semext}  Let $\lambda\colon G\to G'$ be a homomorphism
of connected reductive algebraic groups which maps the connected component
of the centre of $G$ into the connected component of the centre of $G'$.
If $\fE$ is a semistable principal Higgs $G$-bundle, and $\fE'$ is
obtained by extending the structure group $G$ to $G'$ by $\lambda$, then $\fE'$ is semistable.
\end{prop}
\begin{proof} By composing the adjoint representation of $G'$ with the homomorphism
$\lambda$ we obtain a representation $\rho\colon G\to \operatorname{Aut}(\mathfrak g')$; the principal Higgs bundle
obtained by extending the structure group of $\fE$ to $\operatorname{Aut}(\mathfrak g')$ is the bundle of linear
frames of $\Ad(\fE')$ with its natural Higgs field. By Lemma \ref{adjsemi}, this
bundle is semistable, so that $\Ad(\fE')$ is semistable as well. Again by Lemma  \ref{adjsemi},
$\fE'$ is semistable.
\end{proof}

We can now prove a version of Miyaoka's semistability criterion
which works for principal Higgs bundles on projective varieties of
any dimension.

\begin{thm}\label{Miyahigher} Let $\fE$ be a  principal Higgs $G$-bundle $\fE=(E,\phi)$ on $X$. The following conditions are equivalent:
\begin{enumerate} \item
for every parabolic subgroup $P\subset G$ and any nontrivial dominant character $\chi$ of $P$, the line bundle $L_{\chi}^\ast$ restricted to ${\fR}_P(\fE)$ is numerically effective;
\item for every morphism $f\colon C\to X$, where $C$ is a smooth projective curve,  the pullback
 $f^\ast(\fE)$   is semistable;
\item $\fE$ is semistable and $c_2(\Ad(E))=0$ in $H^4(X,\R)$.
\end{enumerate}
\end{thm}
\begin{proof}
Assume that condition (i) holds, and let $f\colon C\to X$ be as in the statement.
The line bundle $L'_\chi$ on $f^\ast(E)/P$ given by the character $\chi$
is a pullback of $L_\chi$. Then $L'_{\chi\vert  {\fR}_P(f^\ast\fE)}$ is nef,
so that by Theorem \ref{princMiya},  $f^\ast(\fE)$  is semistable. Thus (i) implies (ii).

We prove that (ii) implies (iii). Since $f^\ast(\fE)$ is semistable,
by Lemma \ref{adjsemi}, the adjoint Higgs bundle $\Ad(f^\ast(\fE))$ is semistable.
By results proved in \cite{BG1} we have
that $\Ad(\fE)$ is semistable, and
$$\Delta(\Ad(E)) = c_2 (\Ad(E)) - \frac{r-1}{2r}(c_1 (\Ad(E)))^2=0\,.$$
As $G$ is reductive, we
have $\Ad(E)\simeq\Ad(E)^\ast$, so that $c_1 (\Ad(E))=0$, and the previous equation reduces to
$ c_2 (\Ad(E)) =0$. Again using Lemma \ref{adjsemi},
we have that $\fE$ is semistable.

Next we prove that (iii) implies (ii).
This is proved by reversing the previous arguments: $\Ad(\fE)$ is semistable by  Lemma \ref{adjsemi}; thus,  since $c_2(\Ad(E))=0$, by results in \cite{BG1} the Higgs vector bundle
$\Ad(f^\ast(\fE))$ is semistable, and then $f^\ast(\fE)$  is semistable
by Lemma \ref{adjsemi}.

Finally, we show that (ii) implies (i). Let $C'$ be a curve in
${\fR}_P(\fE)$. If it is contained in a fibre of the projection $\pi_P\colon {\fR}_P(\fE)\to X$,
since $\chi$ is dominant, we have $c_1(L_\chi^\ast)\cdot[C']\ge 0$. So we
may assume that $C'$ is not in a fibre. Moreover, 
possibly by replacing it with its normalization,
we may assume it is smooth.   The projection
of $C'$ to $X$ is a finite cover  $\pi_P\colon C'\to C$ to its
image $C$. We may choose a smooth projective curve $C''$ and a
morphism $h\colon C''\to C$ such that $\tilde C=C''\times_CC'$ is
a split unramified cover. Then every sheet $C_j$ of $\tilde C$ is
the image of a section $\sigma_j$ of ${\fR}_P(h^\ast\fE)$. Since
$h^\ast\fE$ is semistable by Lemma \ref{basechange}, we have $\deg
\sigma_j^\ast (L^\ast_\chi)\ge 0$ by Lemma \ref{RamaHiggs}. This
implies (i).
\end{proof}

\begin{corol} Assume that $\fE=(E,\phi)$ is a principal Higgs $G$-bundle,
$\psi\colon G\to G'$ is a surjective group homomorphism, $\fE'=(E',\phi')$ is a principal Higgs $G'$-bundle, and $f\colon E\to E'$ is a $\psi$-equivariant morphism of principal Higgs bundles.
If $\fE$ satisfies one of the conditions of Theorem \ref{Miyahigher}, so does $\fE'$.
\end{corol}
\begin{proof} If $P'$ is a parabolic subgroup of $G'$, then $P'=\psi(P)$ for
a parabolic $P$ in $G$. If $\chi'\colon P'\to \C^\ast$ is a dominant character of $P'$,
the composition $\chi=\chi'\circ\psi$ is a dominant character of $P$. If $f\colon E/P\to E'/P'$
is the induced morphism, we know that $f({\fR}_P(\fE))\subset {\fR}_{P'}(\fE')$, so that
$f^\ast(L^\ast_{\chi'\vert {\fR}_{P'}(\fE')})\simeq L^\ast_{\chi\vert {\fR}_P(\fE)}$. Since
$L^\ast_{\chi\vert {\fR}_P(\fE)}$ is nef, and $f\colon  {\fR}_P(\fE))\to {\fR}_{P'}(\fE')$ is surjective,
$L^\ast_{\chi'\vert {\fR}_{P'}(\fE')}$ is nef as well \cite{Fu}.
\end{proof}

\bigskip\section{Another semistability criterion}\label{another}
In this section we prove another semistability criterion, which
states that a principal Higgs bundle $\fE=(E,\phi)$ is semistable  and $c_2(\Ad(E))=0$ if and only if
the adjoint Higgs bundle $\Ad(\fE)$ is numerically flat, in a sense that we introduced in \cite{BG1}.
This calls for a brief reminder of the notion of numerical effectiveness for Higgs bundles
(a notion that we shall call ``H-nefness''). Let us remark beforehand that
our definition of H-nefness for Higgs bundles
requires to consider Higgs bundles on  singular schemes.
For such spaces there is well-behaved theory of the de Rham complex
\cite{EGAIV}, which is all one needs to define Higgs bundles.

Let $X$ be a scheme over the complex numbers, and $E$  a rank $r$
vector bundle on $X$. For every positive integer $s$ less than
$r$, let $\grass_s(E)$ denote the Grassmann bundle of rank $s$ quotients of
$E$, with projection $p_s : \grass_s(E) \to X$.
There is a universal exact sequence of vector bundles on $\grass_s(E)$
\begin{equation}\label{univeq}
0 \to S_{r-s,E} \xrightarrow{\psi} p_s^*(E) \xrightarrow{\eta}
Q_{s,E} \to 0
\end{equation}
where   $S_{r-s,E}$  is the universal
rank $r-s$ subbundle and $Q_{s,E}$ is the universal rank $s$ quotient
bundle \cite{Ful}.

Given a Higgs bundle $\fE $, we construct closed subschemes
$\hgrass_s(\fE)\subset \grass_s(E)$ pa\-ram\-e\-tr\-iz\-ing rank
$s$ locally-free Higgs quotients (these are the counterparts in
the vector bundle case of the schemes of Higgs reductions we have
introduced previously). We define $\hgrass_s(\fE)$ as the closed
subscheme of $\grass_s(E)$ where the composed morphism
$$(\eta\otimes1)\circ p_s^\ast(\phi) \circ \psi\colon S_{r-s,E}\to
Q_{s,E}\otimes
 p_s^\ast(\Omega_X)$$ vanishes. (This scheme has already been introduced by Simpson in
 a particular case, e.g., semistable Higgs bundles with vanishing Chern classes, cf.~\cite[Cor.~9.3]{Si2}.)\ We denote by $\rho_s$ the projection $\hgrass_s(\fE)\to
X$. The restriction of \eqref{univeq} to the scheme $\hgrass_s(\fE)$
provides the exact sequence of vector bundles
\begin{equation} \label{univg}
0 \to S_{r-s,\fE}\to \rho_s^\ast (\fE) \to Q_{s,\fE}\to 0 \,.
\end{equation}
The Higgs morphism $\phi$ of $\fE$ induces by pullback
a Higgs morphism $\Phi\colon  \rho_s^\ast (\fE) \to
 \rho_s^\ast (\fE) \otimes \Omega_{\hgrass_s(\fE)}$.
Due to the condition
 $(\eta\otimes1)\circ p_s^\ast(\phi) \circ \psi=0$ which is satisfied
 on $\hgrass_s(\fE)$, the morphism $\Phi$ sends $ S_{r-s,\fE}$ to $ S_{r-s,\fE}\otimes \Omega_{\hgrass_s(\fE)}$.  As a result, $S_{r-s,\fE}$ is a Higgs subbundle
 of $\rho_s^\ast (\fE)$, and the quotient $ Q_{s,\fE}$ has a structure of Higgs bundle. Thus \eqref{univg} is   an exact sequence
 of Higgs bundles.

We recall from \cite{BG1} the notion of H-nef Higgs bundle.

\begin{defin}  \label{moddef} A Higgs bundle $\fE$ of rank one is said
to be Higgs-numerically effective (for short, H-nef) if it is
numerically effective in the usual sense. If $\rk \fE \geq 2$ we
require that:
\begin{enumerate} \item all bundles $Q_{s,\fE}$ are Higgs-nef;
\item the line bundle $\det(E)$ is nef.
\end{enumerate}
If both $\fE$ and $\fE^\ast$ are Higgs-numerically effective,
$\fE$ is said to be Higgs-numerically flat (H-nflat).
\end{defin}

We are now in position to state and prove the additional semistability criterion we promised.

\begin{thm} Let $\fE$ be a  principal Higgs bundle $\fE=(E,\phi)$ on a polarized smooth complex projective variety $X$. The following conditions are equivalent.
\begin{enumerate}
\item $\fE$ is semistable and $c_2(\Ad(E))=0$ in $H^4(X,\R)$;
\item the adjoint Higgs bundle $\Ad(\fE)$ is H-nflat.
\end{enumerate}
\end{thm}
\begin{proof} At first we prove this theorem when $X$ is a curve. In this case the claim is the following: $\fE$ is semistable
if and only if $\Ad(\fE)$ is H-nflat. In view of Lemma \ref{adjsemi}, this amounts to proving that $\Ad(\fE)$ is semistable
if and only if it is H-nflat. Since $c_1(\Ad(E))=0$ this holds true (\cite{BG1}, Corollaries 3.4 and 3.6).

Let us assume now that $\dim(X)>1$. If condition (i) holds, then $\fE_{\vert C}$ is semistable for any embedded curve
$C$ (as usual, if $C$ is not smooth one replaces it with its normalization). Thus $\Ad(\fE)_{\vert C}$ is
semistable, hence H-nflat. But this implies that $\Ad(\fE)$ is H-nflat as well.

Conversely, if $\Ad(\fE)$ is H-nflat, then it is semistable (see \cite{BG1}), so that $\fE$ is semistable.
Moreover, all Chern classes of H-nflat Higgs bundles vanish, so that $c_2(\Ad(E))=0$.
\end{proof}

\begin{remark} This characterization shows that the numerically flat principal $G$-bundles defined in \cite{BS} for semisimple structure groups $G$ are no more than the class of principal bundles singled out by one of the equivalent conditions of Theorem \ref{Miyahigher} (cf.~\cite[Thm.~2.5]{BS}).
\end{remark}
\bigskip\section{Some Tannakian considerations}\label{tannaka}
In this section we place our main Theorem \ref{Miyahigher}  into the framework
of Tannakian categories. We recall (see e.g.~\cite{DM})  that a neutral Tannakian category $\cat T$ over
a field $\Bbbk$
is a rigid abelian (associative and commutative) $\Bbbk$-linear tensor category such  that
\begin{enumerate}
\item for every
unit object 1 in $\cat T$, the endomorphism space $\operatorname{End}(1)$
is isomorphic to $\Bbbk$;
\item there is an exact faithful functor $\omega\colon \cat T\to \cat{Vect}_{\Bbbk}$,
called a {\em fibre functor.}
\end{enumerate}
Here $\operatorname{Vect}_{\Bbbk}$ is the category of vector spaces over $\Bbbk$. The standard
example of a neutral Tannakian category is the category $\cat{Rep}(G)_{\Bbbk}$ of  $\Bbbk$-linear representations of
an affine group scheme $G$. Indeed, any neutral Tannakian category can be represented
as $\cat{Rep}(G)_{\Bbbk}$ where $G$ is the automorphism group of the fibre functor $\omega$.
Let $\fE$ be a principal Higgs $G$-bundle
on a (say) complex projective manifold $X$.
For any finite-dimensional linear representation $\rho\colon G\to\operatorname{Aut}(W)$ let
$\mathfrak W = \fE\times_\rho W$ be the associated Higgs vector bundle. This correspondence
defines a \emph{$G$-torsor} on the category $\cat{Higgs}_X$ of Higgs vector bundles on $X$, i.e.,
a faithful and exact functor $\fE\colon\cat{Rep}(G)_{\Bbbk}\to\cat{Higgs}_X$ \cite{Si1}. In general,
this is not always compatible with semistability, i.e., $\fE(\rho,W)$ is not always semistable
even when $\fE$ is. In order to have that, we need to impose some conditions. For instance,
we may assume  that  every representation $\rho\colon G\to\operatorname{Aut}(W)$ maps the
connected component of the centre of $G$ containing the identity to the centre
of $\operatorname{Aut}(W)$ (this happens e.g.~when $G$ is semisimple). When this is true, we
say that $G$ is \emph{central}.

Let $\cat{Higgs}_X^\Delta$ be the full subcategory of $\cat{Higgs}_X$ whose
objects $\mathfrak W$ satisfy $\Delta(\mathfrak W)=0$ and are semistable. Since
$$\Delta (\tilde V\otimes \tilde W) = \rk(\tilde W)\Delta(\tilde V) +  \rk(\tilde V)\Delta(\tilde W)\,,$$
and the tensor product of semistable Higgs bundles is semistable \cite{Si1},
it is a tensor category.  However, it is not additive but only preadditive. Let $\cat{Higgs}_X^{\Delta,+}$ 
be its additive completion (see e.g.~\cite{Freyd}). We may now prove the following characterization. \begin{prop} Assume that $G$ is central. There is one-to-one correspondence between principal Higgs $G$-bundles $\fE$ satisfying one of the conditions
of Theorem \ref{Miyahigher} and $G$-torsors on the category $\cat{Higgs}_X$ taking values
in   $\cat{Higgs}_X^{\Delta,+}$. \end{prop}
\begin{proof} Given a principal Higgs $G$-bundle $\fE$ and a representation  $\rho\colon G\to\operatorname{Aut}(W)$ the associated Higgs vector bundle $\mathfrak W$  is semistable  and satisfies
$\Delta(\mathfrak W)=0$  since $\Delta(\mathfrak W)$ is a multiple of $c_2(\Ad(E))$. Conversely, given a $G$-torsor on $\cat{Higgs}_X^{\Delta,+}$, one builds a principal Higgs $G$-bundle $\fE$ as in \cite[Ch.~6]{Si1}. We prove that $\fE$ is
semistable. If $\mathfrak W$
is an associated Higgs vector  bundle  via a faithful representation, $\Ad(\fE)$ is a Higgs subbundle of
$\End(\mathfrak W)$. If $\mathfrak W$ is semistable,
since  $c_1(\Ad(\fE))=c_1(\End(\mathfrak W))=0$ the bundle   $\Ad(\fE)$ is semistable, so that
 $\fE$ is semistable as well. To prove $c_2(\Ad(E))=0$ it is enough to choose for $\rho$ the adjoint representation.
\end{proof}

\end{document}